\documentstyle[12pt]{article}
\def\Bbb R{{\rm \bf R}}
\def\proclaim#1{\vskip2mm{\bf #1}\em}
\def\endproclaim{\em \vskip2mm}
\def\tag#1{\eqno(#1)}
\def\gathered{\begin{array}{c}}
\def\endgathered{\end{array}}
\def\text{\mbox}
\begin{document}

\title {A note on the enclosure method for an inverse obstacle scattering problem
with a single point source}
\author{Masaru IKEHATA\footnote{
Department of Mathematics,
Graduate School of Engineering,
Gunma University, Kiryu 376-8515, JAPAN}}
\date{   }
\maketitle
\begin{abstract}
This paper gives a note on an application of the enclosure method to
an inverse obstacle scattering problem governed by the Helmholtz
equation in two dimensions. It is shown that
one can uniquely
determine the convex hull of an unknown sound-hard polygonal
obstacle from the trace of the total wave that was exerted by a
single point source onto a known circle surrounding the obstacle
provided the source is sufficiently far from the obstacle. The
result contains a formula that extracts the value of the support
function of the obstacle at a generic direction.
Some other applications to thin obstacles, obstacles in a layered medium
and the far-field equation in the linear sampling method are also included.

\noindent
AMS: 35R30


\end{abstract}

%\tableofcontents

\section{Introduction}

The {\it enclosure method} was introduced in \cite{IE} for inverse
boundary value problems for elliptic equations which are motivated by the possibility of applications to
{\it electrical impedance tomography},
{\it diffraction tomography}, etc.. Therein the observation data are
formulated by using the Dirichlet-to-Neumann map (or
Neumann-to-Dirichlet map) associated with the governing equation
of a `signal' propagating inside the medium. It aims at
extracting information about the location and shape of unknown
discontinuity embedded in a known reference medium that gives an effect on the propagation of the signal,
such as an
obstacle, inclusion, crack, etc. from data observed on the
boundary of the medium.  Now we have many applications of this
method, see, e.g., \cite{IE2, ISIL, IMIT, IO3, IS3}.

In \cite{IE3} it was shown that, in a simplified situation
a {\it single set} of the Dirichlet and
Neumann data gives information about the {\it convex hull} of unknown discontinuity.
It was the starting point of the {\it single measurement
version} of the enclosure method
and we have already many applications, e.g.,
\cite{IE4, IE100, IRE, ITRANS, IHER, ILAY2, IH, IO, IO2}.

This paper is closely related to \cite{IE4}.
Therein we considered an inverse obstacle scattering problem
of acoustic wave in two dimensions.
The problem is to reconstruct a two-dimensional
obstacle from the Cauchy data on a circle surrounding the obstacle
of the total wave field generated by a {\it single} incident plane
wave with a {\it fixed} wave number. Let us make a review of one of
the results therein.

We consider  a {\it polygonal} obstacle denoted by $D$, that
is, $D\subset\Bbb R^2$ takes the form $D_1\cup\cdots\cup D_m$ with
$1\le m<\infty$ where each $D_j$ is open and a polygon; $\overline
D_j\cap\overline D_{j'}=\emptyset$ if $j\not=j'$.

The total wave field $u$ outside obstacle $D$ takes the form $u(x;d,k)=e^{ikx\cdot d}+w(x)$
with $k>0$, $d\in S^1$
and satisfies
$$\begin{array}{c}
\displaystyle
\triangle u+k^2 u=0\,\,\text{in}\,\Bbb R^2\setminus\overline D,\\
\\
\displaystyle
\frac{\partial u}{\partial\nu}=0\,\,\text{on}\,\partial D,\\
\\
\displaystyle
\lim_{r\longrightarrow\infty}\sqrt{r}\left(\frac{\partial w}{\partial r}-ikw\right)=0,
\end{array}
$$
where $r=\vert x\vert$ and $\nu$ denotes the unit outward normal relative to $D$.
The last condition above is called the Sommerfeld radiation condition.

Let $B_R$ be an open disc with radius $R$ centered at a fixed point
satisfying $\overline D\subset B_R$.
We assume that $B_R$ is {\it known}.  Our data are  $u=u(\,\cdot\,;d,k)$ and $\partial u/\partial\nu$ on $\partial B_R$
for a fixed $d$ and $k$,
where $\nu$ is the unit outward normal relative to $B_R$.
Let $\omega$ and $\omega^{\perp}$ be two unit
vectors perpendicular to each other.  We always choose the orientation of $\omega^{\perp}$ and $\omega$ coincides with
that of $\mbox{\boldmath $e$}_1$ and $\mbox{\boldmath $e$}_2$ and thus $\omega^{\perp}$ is unique.

We make use of the special complex exponential solution of the
Helmholtz equation $(\triangle+k^2)v=0$ in $\Bbb R^2$:
$$\displaystyle
v_{\tau}(x;\omega)=e^{x\cdot(\tau\omega+i\sqrt{\tau^2+k^2}\omega^{\perp})},\,\,x\in\Bbb R^2,
$$
where $\tau>0$
is a parameter.

Recall the {\it support function} of $D$:  $h_D(\omega)=\sup_{x\in D}x\cdot\omega$.
We say that $\omega$ is {\it regular} with respect to $D$ if the set
$\displaystyle
\partial D\cap\{x\in\Bbb R^2\,\vert\,x\cdot\omega=h_D(\omega)\}$
consists of only one point.

Define
$$\displaystyle
I(\tau;\omega, d, k)
=\int_{\partial B_R}\left(\frac{\partial u}{\partial\nu}v_{\tau}
-\frac{\partial v_{\tau}}{\partial\nu}u\right)
dS.
$$

\proclaim{\noindent Theorem 1.1(\cite{IE4}).}
Assume that $\omega$ is regular with respect to $D$.
Then the formula
$$\displaystyle
\lim_{\tau\longrightarrow\infty}\frac{1}{\tau}
\log
\left\vert I(\tau;\omega, d, k)
\right\vert
=h_D(\omega),
$$
is valid.  Moreover, we have the following:

if $t\ge h_D(\omega)$, then
$\displaystyle
\lim_{\tau\longrightarrow\infty}
e^{-\tau t}\left\vert
I(\tau;\omega, d, k)
\right\vert=0$;

if $t<h_D(\omega)$, then
$\displaystyle
\lim_{\tau\longrightarrow\infty}
e^{-\tau t}\left\vert
I(\tau;\omega, d, k)
\right\vert
=\infty.
$
\endproclaim

In \cite{IHER} a similar formula has been established by using
the {\it far-field
pattern} $F_D(\varphi,d;k)$, $\varphi\in S^1$ of scattered wave $w=u-e^{ikx\cdot d}$ for fixed $d$ and $k$
which determines the leading term of the asymptotic expansion of $w$ as $r\longrightarrow\infty$
in the following sense:
$$\displaystyle
w(r\varphi)\sim\frac{\displaystyle
e^{ikr}}{\displaystyle
\sqrt{r}}F_D(\varphi,d;k).
$$
Moreover, therein instead of volumetric obstacle, similar formulae for {\it thin} sound-hard
obstacle (or {\it screen} )
have also been established with {\it two} incident plane waves.

In this section, we describe another inverse obstacle scattering
problem in which a point source located within a {\it finite distance}
from an unknown
obstacle generates  a scattered wave and one measures the total
wave on a known circle surrounding an unknown obstacle.  One can
see this type of problem in, e.g., a mathematical formulation of
microwave tomography \cite{ S}, subsurface radar \cite{DGS}, etc..

Let $y\in\Bbb R^2\setminus\overline D$.   Let $E=E_D(x,y)$ be the unique solution of the scattering problem:
$$\begin{array}{c}
\displaystyle
(\triangle+k^2)E=0\,\,\text{in}\,\Bbb R^2\setminus\overline D,\\
\\
\displaystyle
\frac{\partial}{\partial\nu}E=-\frac{\partial}{\partial\nu}\Phi_0(\,\cdot\,,y)
\,\text{on}\,\partial D,\\
\\
\displaystyle
\lim_{r\longrightarrow\infty}\sqrt{r}\left(\frac{\partial E}{\partial r}-ikE\right)=0,
\end{array}
$$
where
$$\displaystyle
\Phi_0(x,y)=\frac{i}{4}H^{(1)}_0(k\vert x-y\vert)
$$
and $H^{(1)}_0$ denotes the Hankel function of the first kind \cite{O}.

The total wave outside $D$ exerted by the point source located at $y$ is given by the formula:
$$\displaystyle
\Phi_D(x,y)=\Phi_0(x,y)+E_D(x,y),\,x\in\Bbb R^2\setminus\overline D.
$$

\noindent
{\bf Inverse Problem.}  Let $R_1>R$.
Fix $k>0$ and $y\in\partial B_{R_1}$.
Extract information about the location and shape
of $D$ from $\Phi_D(x,y)$ given at all $x\in\partial B_{R}$.

The aim of this paper is to show that the single measurement
version of the enclosure method still works for this problem.

Define
$$\displaystyle
J(\tau; \omega, y, k)
=
\int_{\partial B_{R}}
\left(
\frac{\partial}{\partial\nu}\Phi_D(x,y)\cdot v_{\tau}(x;\omega)
-\frac{\partial}{\partial\nu}v_{\tau}(x;\omega)\cdot\Phi_D(x,y)\right)dS(x).
$$

The first result of this paper is as follows.

\proclaim{\noindent Theorem 1.2.}
Assume that $\omega$ is regular with respect to $D$ and that
$$\displaystyle
\text{diam}\,D<\text{dist}\,(D,\partial B_{R_1}).
\tag {1.1}
$$
It holds that
$$\displaystyle
\lim_{\tau\longrightarrow\infty}
\frac{1}{\tau}\log
\left\vert J(\tau; \omega, y, k)
\right\vert
=h_D(\omega).
$$
Moreover, we have the following:

if $t\ge h_D(\omega)$, then
$\displaystyle
\lim_{\tau\longrightarrow\infty}e^{-\tau t}\left\vert
J(\tau; \omega, y, k)
\right\vert=0$;

if $t<h_D(\omega)$, then
$\displaystyle
\lim_{\tau\longrightarrow\infty}e^{-\tau t}\left\vert
J(\tau; \omega, y, k)
\right\vert=\infty$.

\endproclaim

It should be pointed out that $(\partial/\partial\nu)\Phi_D(x,y)$ for $x\in\partial B_{R}$ can be computed
from $\Phi_D(x,y)$ for $x\in\partial B_{R}$ by solving the exterior Dirichlet problem for the Helmholtz equation:
$$\begin{array}{c}
\displaystyle
(\triangle+k^2)\tilde{E}=0\,\,\text{in}\,\Bbb R^2\setminus\overline{B_{R}},\\
\\
\displaystyle
\tilde{E}=\Phi_D(\,\cdot\,,y)-\Phi_0(\,\cdot\,,y)\,\,\text{on}\,\partial B_{R},\\
\\
\displaystyle
\lim_{r\longrightarrow\infty}\sqrt{r}\left(\frac{\partial \tilde{E}}{\partial r}-ik\tilde{E}\right)=0.
\end{array}
$$
The computation formula is
$$\displaystyle
\frac{\partial}{\partial\nu}\Phi_D(x,y)=\frac{\partial}{\partial\nu}\Phi_0(x,y)+
\frac{\partial}{\partial\nu}\tilde{E}(x)\,\,\text{on}\,\partial B_{R}.
$$
Needless to say, this kind of remak works also for $\partial u/\partial\nu$ in Theorem 1.1.

Condition (1.1) can be satisfied if $R_1$ is sufficiently
large compared with $R$. It is not known whether condition (1.1)
can be dropped completely. To suggest a possibility next we
present a partial result which does not employ (1.1).

For the description of the second result we introduce special scattered and total fields.
Given $d\in S^1$ choose $\vartheta\in S^1$ in such a way that $\vartheta^{\perp}=d$.
Let $x_0\in\partial D$ and $w=w(x;-d,k,x_0)$ be the unique solution of the scattering problem:
$$\begin{array}{c}
\displaystyle
(\triangle+k^2)w=0\,\,\text{in}\,\Bbb R^2\setminus\overline D,\\
\\
\displaystyle
\frac{\partial w}{\partial\nu}=-\frac{\partial}{\partial\nu}\{(x_0-x)\cdot\vartheta\,e^{-ikx\cdot d}\}\,\,\text{on}\,\partial D,\\
\\
\displaystyle
\lim_{r\longrightarrow\infty}\,\sqrt{r}\left(\frac{\partial w}{\partial r}-ik\,w\right)=0.
\end{array}
\tag {1.2}
$$
Define
$$
\displaystyle
u(x;-d,k,x_0)=(x_0-x)\cdot\vartheta\,e^{-ikx\cdot d}+w(x;-d,k,x_0).
\tag {1.3}
$$
Note that the function $x\longmapsto (x_0-x)\cdot\vartheta e^{-ikx\cdot d}$ satisfies
the Helmholtz equation in the whole plane and
the radiation condition  for $w(x;-d,k,x_0)$ yields that
$$
\displaystyle
u(x;-d,k,x_0)=(x_0-x)\cdot\vartheta e^{-ikx\cdot d}+O(r^{-1/2})
\tag {1.4}
$$
as $r\longrightarrow\infty$.
The second result of this paper is as follows.

\proclaim{\noindent Theorem 1.3.}
Let $\omega$ be regular with respect to $D$
and $x_0\in\partial D$ be the point satisfying $x_0\cdot\omega=
h_D(\omega)$.
Assume that $y\in\partial B_{R_1}$
satisfies $u(y;-d_1,k,x_0)\not=0$ or $u(y;-d_2,k,x_0)\not=0$,
where
$d_1, d_2\in S^1$ are directions of two sides of $D$ that meet at $x_0$.
Then all the conclusions in Theorem 1.2 are valid.
\endproclaim

Note that in this theorem (1.1) is not assumed at the price of
introducing another implicit restriction on the location of $y$
relative to $D$.
From (1.4) we see that given $\delta>0$ and $\vartheta\in S^1$ if $R_0>0$ is sufficiently large,
for all $R_1\ge R_0$ and $y\in\partial B_{R_1}$ with $\vert (x_0-y)\cdot\vartheta\vert\ge\delta$ it holds that
$u(y;-d,k,x_0)\not=0$ with $d=\vartheta^{\perp}$.
We leave the problem of removing the condition completely as an open problem.  See also Remark 2.2 in Section 2.

Since the set of all $\omega$ which is not regular with respect to
given $D$ is finite and the support function of $D$ is continuous,
as a corollary of Theorem 1.2 we have uniqueness of determining
the convex hull of $D$ from $\Phi_D(x,y)$ given at all
$x\in\partial B_{R}$ for a single $y\in\partial B_{R_1}$ and
$k>0$ provided (1.1). It seems that this type of uniqueness with a
single point source had not appeared in the previous study. See
\cite{EY} and references therein for uniqueness of a polygonal
obstacle with the far-field pattern of the scattered wave exerted
by a single plane wave.  It seems that their argument heavily
depends on the fact that the total wave approaches the plane wave
at infinity.  This is not true for the total wave
exerted by a point source
since we have, as $r\longrightarrow\infty$,
$$\displaystyle
\Phi_D(rd,y)
\sim \frac{e^{i\pi/4}}{\sqrt{8\pi k}}\frac{e^{ikr}}{\sqrt{r}}u(y;-d,k).
\tag {1.5}
$$
See \cite{R, IPR, P} for the derivation for $D$ with a smooth boundary.
In our case $\partial D$ is not smooth; however, a minor modification of the
proof still works.
Equation (1.5) means that the far-field pattern of $\Phi_D(x,y)$ as a function of $x$ is given by
$u(y;-d,k)$ multiplied by a known constant.  This formula has been used in the {\it probe method} \cite{IPR}
and the {\it singular sources method} \cite{P1}.

Theorem 1.2 together with (1.5) yields

\proclaim{\noindent Corollary 1.1.}  Let $R_1>R$ and $\overline D\subset B_{R}$.  Assume that $D$ satisfies (1.1).
Fix $k>0$ and $y\in\partial B_{R_1}$.
One can uniquely determine the convex hull of $D$ from the data
$u(y;-d,k)$ given at all $d\in S^1$.

\endproclaim

{\it\noindent Proof.}  The proof is divided into three steps.

\noindent
(i) Use (1.5) to compute the far-field pattern of $\Phi_D(x,y)$ from $u(y;-d,k)$ given at all $d\in S^1$.

\noindent
(ii) Use, e.g., the {\it point source method} \cite{P} to compute $\Phi_D(x,y)$ together with its normal derivative
for $x\in\partial B_{R}$ from the far-field pattern of $\Phi_D(x,y)$.

\noindent
(iii) Use Theorem 1.2 to compute $h_D(\omega)$ for a generic $\omega$ from
$\Phi_D(x,y)$ together with its normal derivative
for $x\in\partial B_{R}$.

\noindent
$\Box$

Summing up, we obtained two procedures for estimating the convex hull of an unknown sound-hard
polygonal obstacle by using two types of the data.

The first type of the data is given by the following process:

(A)  produce the total wave by a {\it fixed} point source located
outside a known circle surrounding an unknown obstacle and observe
the wave at {\it all points} on the circle.

The second is as follows:

(B)  produce the total waves by incident plane waves for {\it all directions}
and observe the waves at a {\it fixed} point outside a known circle surrounding an unknown obstacle.

Note that these are different from the {\it reciprocity principle} (\cite{CK}) which is the identity
$F_D(\varphi,d;k)=F_D(-d,-\varphi;k)$ since in this identity the incident wave is always a {\it plane
wave} and one observes the scattered wave at {\it infinity}.

A brief outline of this paper is as follows. Theorems 1.2 and 1.3
are proved in Section 2. Both proofs have a common starting point
with the proof of Theorem 1.1 which we recall before describing
subsections 2.1 and 2.2. In those subsections we complete the
proof of Theorems 1.2 and 1.3. In the last section three other
applications are given.  Two of them are concerned with some
extensions to thin obstacles and obstacles in a layered medium.
In the last of the applications we consider the far-field equation
which plays the central role in the {\it linear sampling method}
\cite{CK1}.  We show that a modification of the argument for the
proof of Theorem 1.2 gives {\it unsolvabilty} of the far-field
equation for polygonal obstacles.

\section{Proof of Theorems 1.2 and 1.3}

First we follow the argument for the proof of Theorem 1.1 (see also \cite{ITRANS}).
For simplicity of notation we set $u(x)=\Phi_D(x,y)$.

Let $x_0$ denote the single point of the set
$\{x\,\vert\,x\cdot\omega=h_D(\omega)\}\cap\partial D$. $x_0$ has
to be a vertex of $D_j$ for some $j$. In what follows we denote by
$B_R(x_0)$ the open disc with radius $R$ centered at $x_0$. Let
$\Theta$ denote the outside angle of $D$ at $x_0$.
$\Theta$ satisfies $\pi<\Theta<2\pi$ since $\omega$ is
regular with respect to $D$.

If one chooses a
sufficiently small $\eta>0$, then one can write
$$\begin{array}{c}
\displaystyle
B_{2\eta}(x_0)\cap (B_{R_1}\setminus\overline D)
=\{x_0+r(\cos\,\theta\text{\boldmath $a$}+\sin\,\theta\text{\boldmath $a$}^{\perp})\,\vert\,
0<r<2\eta,\,0<\theta<\Theta\},
\\
\\
\displaystyle
B_{\eta}(x_0)\cap\partial D=\Gamma_p\cup \Gamma_q\cup\{x_0\}
\end{array}
$$
where $\displaystyle
\text{\boldmath $a$}=\cos\,p\,\omega^{\perp}+\sin\,p\,\omega$,
$\displaystyle
\text{\boldmath $a$}^{\perp}=-\sin\,p\,\omega^{\perp}+\cos\,p\,\omega$;
$\displaystyle -\pi<p<0$;
$\displaystyle
\Gamma_p=\{x_0+r\text{\boldmath $a$}\,\vert\,0<r<\eta\}$,
$\displaystyle\Gamma_q=\{x_0+r(\cos\,\Theta\text{\boldmath $a$}+\sin\,\Theta\text{\boldmath $a$}^{\perp})\,\vert\,0<r<\eta\}$.
Note that the orientation of $\text{\boldmath $a$}$, $\text{\boldmath $a$}^{\perp}$ coincides with that of
$\text{\boldmath $e$}_1, \text{\boldmath $e$}_2$.
See also Figure 1 of \cite{IE3}.

The quantity $-p$ means the angle between two vectors $\omega^{\perp}$ and $\text{\boldmath $a$}$.
$p$ satisfies $\Theta>\pi+(-p)$. Set $q=\Theta-2\pi+p$.  Then we have $-\pi<q<p<0$ and
the expression
$$\begin{array}{c}
\displaystyle
\Gamma_p=\{x_0+r(\cos\,p\,\omega^{\perp}+\sin\,p\,\omega)\,\vert\,0<r<\eta\},\\
\\
\displaystyle
\Gamma_q=\{x_0+r(\cos\,q\,\omega^{\perp}+\sin\,q\,\omega)\,\vert\,0<r<\eta\}.
\end{array}
$$
This is the meaning of $p$ and $q$.

We set
$$\displaystyle
u(r,\theta)=u(x),\,\,x=x_0+r(\cos\,\theta\,\text{\boldmath $a$}+\sin\,\theta\,
\text{\boldmath $a$}^{\perp}).
$$

The $u$ can be expanded as
$$\displaystyle
u(r,\theta)=\alpha_1J_0(kr)+\sum_{n=2}^{\infty}\alpha_n J_{\lambda_n}(kr)\cos \lambda_n\theta,\,\,0<r<\eta, 0<\theta<\Theta,
$$
where the $\lambda_n$ describes the {\it singularity} of $u$ as $r\longrightarrow 0$ and in this case explicitly given by
the formula $\displaystyle\lambda_n=(n-1)\pi/\Theta$, $J_{\lambda_n}$ stands for the Bessel function of order $\lambda_n$.

One of key points is {\it
introducing a new parameter} $s$ instead of $\tau$ by the equation
$\displaystyle s=\sqrt{\tau^2+k^2}+\tau$, we obtain, as
$s\longrightarrow\infty$, the {\it complete} asymptotic expansion
$$\displaystyle
J(\tau;\omega,y,k)
\,e^{-i\sqrt{\tau^2+k^2}x_0\cdot\omega^{\perp}}e^{-\tau h_D(\omega)}
\sim
-i\sum_{n=2}^{\infty}\frac{e^{i\frac{\pi}{2}\lambda_n}k^{\lambda_n}\alpha_nK_n}
{s^{\lambda_n}},
\tag {2.1}
$$
where $K_n$ are constants given by the formula
$$\displaystyle
K_n=e^{ip\lambda_n}+(-1)^n e^{iq\lambda_n}.
$$
For the derivation of this expansion see \cite{IE4}.  Note that constants $K_n$ are exactly same
as the corresponding ones
in \cite{IE3, IE4}.

Now all the statements in Theorem 1.2 follow from (2.1) and
another key point: $\exists n\ge 2$\,\, $\alpha_nK_n\not=0$.
This is due to a contradiction argument.
Assume that the assertion is not true, that is,
$\forall n\ge 2$\,\,$\alpha_nK_n=0$.

\noindent
{\bf Case A.}
First we consider the case when $\Theta/\pi$ is {\it irrational}.
It is easy to see that $K_n\not=0$ for all $n\ge 2$.  Thus,
$\alpha_n=0$ and this yields $u(r,\theta)=\alpha_1J_0(kr)$
for $0<r<\eta$ and $0<\theta<\Theta$.
Since this right-hand side is an entire solution of the
Helmholtz equation, the unique continuation property of the
solution of the Helmholtz equation yields $u(x)=\alpha_1J_0(k\vert
x-x_0\vert)$ in $\Bbb R^2\setminus\overline D$.
This implies that $u$ has to be bounded in a neighbourhood of $y$.
However, since $$\displaystyle
\Phi_0(x,y)\sim\frac{1}{2\pi}\log\frac{1}{\vert x-y\vert}
$$
as $x\longrightarrow y$ and $E_D(\,\cdot\,,y)$ is smooth in a neighbourhood of $y$,
one knows that $u(x)=\Phi_0(x,y)+E_D(x,y)$ is not bounded in any neighbouhood of $y$.
Contradiction.

\noindent
{\bf Case B.}
Next consider the case when $\Theta/\pi$ is {\it rational}.
One can write
$$\displaystyle
\frac{\Theta}{\pi}=1+\frac{b}{a},
$$
where $a(\ge 2)$ and $b(\ge 1)$ are integers and {\it mutually prime}.
Then we have
$$\displaystyle
\{n\ge 2\,\vert\,K_n=0\}=\{1+l(a+b)\,\vert\,l=1,2,\cdots\}.
\tag {2.2}
$$
Note also that $\lambda_{1+l(a+b)}=al$, $l=1,2,\cdots$.
From the assumption of the contradiction argument one knows if $n$
satisfies $K_n\not=0$, then $\alpha_n=0$.  From this together with (2.2)
we have
$$\displaystyle
u(r,\theta)=\sum_{l=0}^{\infty}
\alpha_{1+l(a+b)}J_{al}(kr)\cos\,al\theta.
\tag {2.3}
$$

Hereafter we take {\it two courses} corresponding to Theorems 1.2
and 1.3.

\subsection{Completion of the proof of Theorem 1.2}

Since each $al$ in (2.3) is an integer, the right-hand side of (2.3) gives
a continuation of $u$ onto $(\Bbb R^2\setminus(\overline D\cup\{y\}))\cup B_{\eta}(x_0)$
as a solution of the Helmholtz equation and the continuation which we denote by $\tilde{u}$
satisfies the rotation invariance in
$B_{\eta}(x_0)$:
$$\displaystyle
\tilde{u}\left(r,\theta+\frac{2\pi}{a}\right)
=\tilde{u}(r,\theta).
\tag {2.4}
$$

Now having (1.1) and (2.4),  one can apply Friedman-Isakov's
extension argument \cite{FI} to $\tilde{u}$.  See also
\cite{ITRANS} for the detail of the argument applied to a
penetrable obstacle case. As a result one gets a continuation of
$\tilde{u}$ onto $\Bbb R^2\setminus\{y\}$ as a solution of the
Helmholtz equation.  Since
$\tilde{u}(x)=u(x)=\Phi_0(x,y)+E_D(x,y)$ in $\Bbb
R^2\setminus\overline D$ and $E_D(x,y)$ is smooth in a
neighbourhood of $y$, one concludes that $E_D(\,\cdot\,,y)$ can be
continued as a solution of the Helmholtz equation in $\Bbb R^2$.
The continuation satisfies the Sommerfeld radiation condition and
therefore has to be identically zero. This gives, in particular,
$\tilde{u}(x)=\Phi_0(x,y)$ in $B_{\eta}(x_0)$ and from (2.4) one
gets
$$\displaystyle
\Phi_0(x_0+rz(\theta),y)
=\Phi_0\left(x_0+rz\left(\theta+\frac{2\pi}{a}\right),y\right),\,\,0<r<\eta,\,\,\theta\in\Bbb R,
\tag {2.5}
$$
where
$\displaystyle
z(\theta)=\cos\,\theta\,\mbox{\boldmath $a$}+\sin\,\theta\,\mbox{\boldmath $a$}^{\perp}$.

Since both sides of (2.5) satisfy the same Helmholtz equation in
$\vert x-x_0\vert<\vert y-x_0\vert$, the unique continuation
property of the solution of the Helmholtz equation yields:  (2.5)
is valid for all $r$ with $r<\vert y-x_0\vert$. Choose a
$\theta_0$ in such a way that $y=x_0+\vert
y-x_0\vert\,z(\theta_0)$. Since $2\pi/a\le\pi$, we have
$y\not=x_0+\vert y-x_0\vert\,z(\theta_0+2\pi/a)$.  Then letting
$\theta=\theta_0$ and $r\uparrow \vert y-x_0\vert$ in (2.5), we
have a contradiction because of the singularity of $\Phi_0(x,y)$ as
$x\longrightarrow y$.

This completes the proof of Theorem 1.2.

\noindent
$\Box$

{\bf\noindent Remark 2.1.}
In the proof of Theorem 1.1, we never
make use of (2.4) after having (2.3) and instead take another
course for the total field $u=u(\,\cdot\,;d,k)$ in Theorem 1.1.

The argument is as follows.
Set
$$\begin{array}{c}
\displaystyle
d_1=\cos\,\theta\,\text{\boldmath $a$}+\sin\,\theta\,\text{\boldmath $a$}^{\perp}\vert_{\theta=\Theta-\pi},\,\,
\vartheta_1=\cos\,\theta\,\text{\boldmath $a$}+\sin\,\theta\,\text{\boldmath $a$}^{\perp}\vert_{\theta=(\Theta-\pi)+\pi/2},
\\
\\
\displaystyle
d_2=\cos\,\theta\text{\boldmath $a$}+\sin\,\theta\text{\boldmath $a$}^{\perp}\vert_{\theta=\pi},\,\,
\vartheta_2=\cos\,\theta\,\text{\boldmath $a$}+\sin\,\theta\,\text{\boldmath $a$}^{\perp}\vert_{\theta=\pi+\pi/2}.
\end{array}
\tag {2.6}
$$
Note that $d_1$ and $d_2$ are directed along the two sides that meet at $x_0$.

From the right-hand
side of (2.3) one gets: for all $r$ with $0<r<<1$ $\displaystyle
\nabla u(x_0+rd_1)\cdot\vartheta_1=0$ and
$\displaystyle\nabla u(x_0+rd_2)\cdot\vartheta_2=0$. Then a reflection argument
in \cite{AD} yields that this is true for all $r>0$.  However, from
this together with the asymptotic behaviour of $\nabla u\sim\nabla e^{ikx\cdot d}$
as $r\longrightarrow\infty$ one gets $d\cdot\vartheta_1=d\cdot\vartheta_2=0$.
Contradiction.

The advantage of this argument is: one does not need to use (1.1).
In the following subsection we employ this argument after (2.3).

\subsection{Completion of the proof of Theorem 1.3}

We use the same notation as (2.6).
First we claim that, as $r\longrightarrow\infty$,
$$\displaystyle
\nabla_x\Phi_D(x_j,y)\cdot\vartheta_j
=\frac{ik^{1/2}}{4}\sqrt{\frac{2}{\pi}}e^{i\pi/4}e^{ik(x_0-y)\cdot d_j}
\frac{e^{ikr}}{r^{3/2}}
u(y;-d_j,k,x_0)+O\left(\frac{1}{r^{5/2}}\right),
\tag {2.7}
$$
where $x_j=x_0+rd_j$.

This is proved as follows. Total field $\Phi_D(\,\cdot\,,y)$ has
the expression
$$\displaystyle
\Phi_D(x,y)
=\Phi_0(x,y)
+\int_{\partial D}\frac{\partial}{\partial\nu(z)}\Phi_0(z,x)\Phi_D(z,y)dS(z),
\,\,x\in\Bbb R^2\setminus\overline D.
\tag {2.8}
$$
Since
$$\displaystyle
\nabla_x\Phi_0(x,y)
=\frac{ik}{4}\frac{x-y}{\vert x-y\vert}(H^{(1)}_0)'(k\vert x-y\vert),
$$
we have
$$\displaystyle
\nabla_x\Phi_0(x_j,y)\cdot\vartheta_j
=\frac{ik}{4}\frac{(x_j-y)\cdot\vartheta_j}{\vert x_j-y\vert}
(H^{(1)}_0)'(k\vert x_j-y\vert).
$$
Here we note that
$(x_j-y)\cdot\vartheta_j=(x_0-y)\cdot\vartheta_j+rd_j\cdot\vartheta_j
=(x_0-y)\cdot\vartheta_j$ since $d_j\cdot\vartheta_j=0$.
This gives
$$\displaystyle
\nabla_x\Phi_0(x_j,y)\cdot\vartheta_j
=\frac{ik}{4}\frac{(x_0-y)\cdot\vartheta_j}{\vert x_j-y\vert}
(H^{(1)}_0)'(k\vert x_j-y\vert).
$$

By (4.03) on p.238 in \cite{O}, we know that as $r\longrightarrow\infty$,
$H^{(1)}_0(r)$ and its derivatives satisfy
$$\displaystyle
H^{(1)}_0(r)=\sqrt{\frac{2}{\pi}}e^{-i\pi/4}r^{-1/2}e^{ir}
\left(1+O(r^{-1})\right),
$$
$$\displaystyle
(H^{(1)}_0)'(r)
=\sqrt{\frac{2}{\pi}}e^{i\pi/4}r^{-1/2}e^{ir}
\left(1+O(r^{-1})\right),
\tag {2.9}
$$
and
$$\displaystyle
(H^{(1)}_0)''(r)
=i\sqrt{\frac{2}{\pi}}e^{i\pi/4}r^{-1/2}e^{ir}
\left(1+O(r^{-1})\right).
\tag {2.10}
$$

Using
$$\displaystyle
\vert x_j-y\vert
=r+(x_0-y)\cdot d_j+O(r^{-1}),
$$
and (2.9), we have
$$\displaystyle
(H^{(1)}_0)'(k\vert x_j-y\vert)
=
\sqrt{\frac{2}{\pi}}e^{i\pi/4}k^{-1/2}r^{-1/2}e^{ikr}
e^{ik(x_0-y)\cdot d_j}
\left(1+O(r^{-1})\right)
$$
and thus
$$\begin{array}{c}
\displaystyle
\nabla_x\Phi_0(x_j,y)\cdot\vartheta_j
=
\frac{ik^{1/2}}{4}\frac{(x_0-y)\cdot\vartheta_j}{r^{3/2}}
\sqrt{\frac{2}{\pi}}e^{i\pi/4}e^{ikr}
e^{ik(x_0-y)\cdot d_j}+O(r^{-5/2}).
\end{array}
\tag {2.11}
$$

Let $z\in\partial D$.  We have
$$\displaystyle
\frac{\partial}{\partial\nu(z)}
\Phi_0(z,x)
=-\frac{ik}{4}\frac{(x-z)\cdot\nu(z)}{\vert x-z\vert}
(H^{(1)}_0)'(k\vert x-z\vert)
$$
and thus
$$\begin{array}{c}
\displaystyle
\nabla_x\left(\frac{\partial}{\partial\nu(z)}
\Phi_0(z,x)\right)
=-\frac{ik}{4}\frac{\nu(z)}{\vert x-z\vert}
(H^{(1)}_0)'(k\vert x-z\vert)\\
\\
\displaystyle
+\frac{ik}{4}
\frac{(x-z)\cdot\nu(z)\,(x-z)}
{\vert x-z\vert^3}
(H^{(1)}_0)'(k\vert x-y\vert)
-\frac{ik^2}{4}\frac{(x-z)\cdot\nu(z)\,(x-z)}{\vert x-z\vert^2}
(H^{(1)}_0)''(k\vert x-z\vert).
\end{array}
$$
This together with $(x-z)\cdot\vartheta_j=(x_0-z)\cdot\vartheta_j$ yields
$$\begin{array}{c}
\displaystyle
\nabla_x\left(\frac{\partial}{\partial\nu(z)}
\Phi_0(z,x)\right)\vert_{x=x_j}\cdot\vartheta_j
=-\frac{ik}{4}\frac{\nu(z)\cdot\vartheta_j}{\vert x_j-z\vert}
(H^{(1)}_0)'(k\vert x_j-z\vert)\\
\\
\displaystyle
+\frac{ik}{4}
\frac{(x_j-z)\cdot\nu(z)\,(x_0-z)\cdot\vartheta_j}
{\vert x_j-z\vert^3}
(H^{(1)}_0)'(k\vert x_j-z\vert)\\
\\
\displaystyle
-\frac{ik^2}{4}\frac{(x_j-z)\cdot\nu(z)\,(x_0-z)\cdot\vartheta_j}{\vert x_j-z\vert^2}
(H^{(1)}_0)''(k\vert x_j-z\vert).
\end{array}
$$
Note that the second term of this right-hand side is estimated by $O(r^{-5/2})$.
It follows from these and (2.10) that
$$\begin{array}{c}
\displaystyle
\nabla_x\left(\frac{\partial}{\partial\nu(z)}
\Phi_0(z,x)\right)\vert_{x=x_j}\cdot\vartheta_j
\\
\\
\displaystyle
=-\frac{ik^{1/2}}{4}
\frac{\nu(z)\cdot\vartheta_j+ikd_j\cdot\nu(z)\,(x_0-z)\cdot\vartheta_j}{r^{3/2}}
\sqrt{\frac{2}{\pi}}e^{i\pi/4}e^{ikr}
e^{ik(x_0-z)\cdot d_j}
+O(r^{-5/2}).
\end{array}
\tag {2.12}
$$
Now from (2.7), (2.8), (2.11) and (2.12), we obtain
$$\displaystyle
\nabla_x\Phi_D(x_j,y)\cdot\vartheta_j
=\frac{ik^{1/2}}{4}\sqrt{\frac{2}{\pi}}e^{i\pi/4}e^{ik(x_0-y)\cdot d_j}
\frac{e^{ikr}}{r^{3/2}}
U(y;d,k,x_0)+O(\frac{1}{r^{5/2}}),
$$
where
$$\begin{array}{c}
\displaystyle
U(y;d,k,x_0)
=(x_0-y)\cdot\vartheta_j e^{ik(x_0-y)\cdot d_j}\\
\\
\displaystyle
-\int_{\partial D}
\left(\nu(z)\cdot\vartheta_j+ik d_j\cdot\nu(z)\,(x_0-z)\cdot\vartheta_j\right)
e^{ik(x_0-z)\cdot d_j}\Phi_D(z,y)dS(z).
\end{array}
\tag {2.13}
$$
Define
$$\displaystyle
\Psi_j(x)=(x_0-x)\cdot\vartheta_j\,e^{-ikx\cdot d_j}.
$$
Since
$$\displaystyle
\frac{\partial}{\partial\nu(z)}
\left((x_0-z)\cdot\vartheta_j\,e^{-ik z\cdot d_j}\right)
=-\left(\nu(z)\cdot\vartheta_j+ikd_j\cdot\nu(z)\,(x_0-z)\cdot\vartheta_j\right)
e^{-ikz\cdot d_ju},
$$
one can rewrite (2.13) as
$$\displaystyle
U(y;d,k,x_0)=\Psi_j(y)
+\int_{\partial D}
\frac{\partial}{\partial\nu(z)}\Psi_j(z)\cdot\Phi_D(z,y)dS(z).
\tag {2.14}
$$
On the other hand, a combination of  Green's identity, the
Sommerfeld radiation condition for $w_j(z)\equiv w(z;-d_j,k,x_0)$
and $\Phi_D(z,y)$ and the boundary condition in (1.2) gives
$$\begin{array}{c}
\displaystyle
w_j(y)=\int_{\partial D}\left(w_j(z)\frac{\partial}{\partial\nu(z)}\Phi_D(z,y)
-\Phi_D(z,y)\frac{\partial}{\partial\nu(z)}w_j(z)\right)dS(z)\\
\\
\displaystyle
=-\int_{\partial D}
\Phi_D(z,y)\frac{\partial}{\partial\nu(z)}w_j(z)dS(z)\\
\\
\displaystyle
=\int_{\partial D}\frac{\partial}{\partial\nu(z)}\{(x_0-z)\cdot\vartheta_j\,e^{-ikz\cdot d_j}\}\Phi_D(z,y)dS(z).
\end{array}
$$

Therefore we see that the left-hand side of (2.14) coincides with
$u(y;-d_j,k,x_0)$.  This completes the proof of (2.7).

Now the proof of Theorem 1.3 starts with having (2.3).
By the same reason described in Remark 2.1, from (2.3) we have
$$\displaystyle
\nabla_x\Phi_D(x_j,y)\cdot\vartheta_j=0,
\tag {2.15}
$$
where $x_j=x_0+rd_j$ and $0<r<<1$.
First consider the case when $y\not=x_0+rd_j$ for all $r>0$.
In this case a reflection argument in \cite{AD} ensures that (2.15) is valid
for all $r>0$ and thus (2.7) yields
$$\displaystyle
u(y;-d_j,k,x_0)=0.
\tag {2.16}
$$

If $y=x_0+\vert y-x_0\vert d_j$, then $\nabla_x\Phi_0(x,y)\cdot\vartheta_j=0$
for $x=x_0+rd_j$ with $0<r<\vert y-x_0\vert$ and $r>\vert y-x_0\vert$.  Then form (2.15) we have
$\nabla_x E_D(x,y)\cdot\vartheta_j=0$ for $x=x_0+rd_j$ with $0<r<<1$ and this is true for
all $r>0$ by a reflection argument in \cite{AD}.  Therefore we again have (2.15) for all $r>\vert y-x_0\vert$
and thus (2.16) too.

Summing up, in any case we obtain equation (2.16) for $j=1,2$.
This is a contradiction.

This completes the proof of Theorem 1.3.

\noindent
$\Box$

{\bf\noindent Remark 2.2.}
However, (2.16) is coming from only the leading term of the  asymptotic expansion (2.7).
Thus our next problems in this direction are as follows:

(i)  determine the complete asymptotic expansion of
$e^{-ikr}\nabla_x\Phi_D(x_j,y)\cdot\vartheta_j$ as
$r\longrightarrow\infty$:
$$\displaystyle
e^{-ikr}\nabla_x\Phi_D(x_j,y)\cdot\vartheta_j
\sim\sum_{m=0}^{\infty}A_m(y,x_0,d_j,k)r^{-(3/2+m)}.
$$

(ii)  If $A_m(y,x_0,d_j,k)=0$ for all $m$, then what happens?

The main obstruction in this approach is the complexity
of computing the asymptotic expansion in (i) as can be seen in the
proof of Theorem 1.3.

\section{Other applications}

In this last section, instead we give three applications of the argument done in Theorems 1.2 and 1.3.

\subsection{Thin obstacle}

It should be pointed out that the advantage of the assumption in
Theorem 1.3 is that the result can be extended to a {\it thin}
obstacle case.

First we review a result in \cite{IE4} which employs a single plane wave as an incident wave
and corresponds to Theorem 1.1. 

Let $\Sigma$ be the union of finitely many disjoint closed piecewise linear segments
denoted by $\Sigma_1, \Sigma_2,\cdots,\Sigma_m$.
Assume that there exists  a simply connected open set $D$
such that $D$ is a polygon and each $\Sigma_j$ consists of sides of $D$.

We assume that $\overline D\subset B_{R}$ with a $R>0$.
We denote by $\nu$ the unit outward normal on $\partial D$ relative to $B_{R}\setminus\overline D$
and set $\nu^+=\nu$ and $\nu^-=-\nu$ on $\Sigma$.
Given $k>0$ and $d\in S^1$ let $u=u(x), x\in\,\Bbb R^2\setminus\Sigma$,
be the solution of the scattering problem
$$\begin{array}{c}
\displaystyle
(\triangle+k^2)u=0\,\,\text{in}\,\,\Bbb R^2\setminus\Sigma,\\
\\
\displaystyle
\frac{\partial u^{\pm}}{\partial\nu^{\pm}}=0\,\,\text{on}\,\Sigma,\\
\\
\displaystyle
\lim_{r\longrightarrow\infty}\,\sqrt{r}\left(\frac{\partial w}{\partial r}-ik\,w\right)=0,
\end{array}
$$
where $w=u-e^{ikx\cdot d}$, $u^+=u\vert_{\Bbb
R^2\setminus\overline D}$ and $u^-=u\vert_{D}$.  Note that this is
a brief description of the problem and for exact one see
\cite{IE4}. Define
$$\displaystyle
I_{\Sigma}(\tau;\omega,d,k)
=\int_{\partial B_{R}}
\left(\frac{\partial u}{\partial\nu} v_{\tau}-\frac{\partial v_{\tau}}{\partial\nu}u\right)dS.
$$

\noindent

In \cite{IE4} we have established the following result.

\proclaim{\noindent Theorem 3.1(\cite{IE4}).} Let $\omega$ be regular with respect to $\Sigma$.
If every end points of $\Sigma_1,\Sigma_2,\cdots,\Sigma_m$ satisfies
$x\cdot\omega<h_{\Sigma}(\omega)$, then the formula
$$\displaystyle
\lim_{\tau\longrightarrow\infty}\frac{1}{\tau}
\log\vert I_{\Sigma}(\tau;\omega,d,k)\vert=h_{\Sigma}(\omega),
$$
is valid. Moreover, we have the following:

if $t\ge h_{\Sigma}(\omega)$, then $\lim_{\tau\longrightarrow\infty}e^{-\tau t}\vert I_{\Sigma}(\tau;\omega,d,k)\vert=0$;

if $t<h_{\Sigma}(\omega)$, then $\lim_{\tau\longrightarrow\infty}e^{-\tau t}\vert I_{\Sigma}(\tau;\omega,d,k)\vert=\infty$.

If there is an end point $x_0$ of some $\Sigma_j$ such that
$x_0\cdot\omega=h_{\Sigma}(\omega)$, then, for $d$ that is not
perpendicular to $\nu$ on $\Sigma_j$ near the point, the same
conclusions as above are valid.

\endproclaim

Note that $\nu$ on $\Sigma_j\cap B_{\eta}(x_0)$ for sufficiently small $\eta>0$
becomes a constant vector if $x_0$ is an end point of $\Sigma_j$.

Here we present a result in which, instead of a single plane wave we make use of
a single point source as an incident wave.

Let $y\in\Bbb R^2\setminus\overline D$.   Let $E=E_{\Sigma}(x,y)$ be the unique solution of the scattering problem:
$$\begin{array}{c}
\displaystyle
(\triangle+k^2)E=0\,\,\text{in}\,\Bbb R^2\setminus\Sigma,\\
\\
\displaystyle
\frac{\partial}{\partial\nu^{\pm}}E^{\pm}=-\frac{\partial}{\partial\nu^{\pm}}\Phi_0(\,\cdot\,,y)
\,\text{on}\,\Sigma,\\
\\
\displaystyle
\lim_{r\longrightarrow\infty}\sqrt{r}\left(\frac{\partial E}{\partial r}-ikE\right)=0.
\end{array}
$$

The total wave outside $\Sigma$ exerted by the point source located at $y$ is given by the formula:
$$\displaystyle
\Phi_{\Sigma}(x,y)=\Phi_0(x,y)+E_{\Sigma}(x,y),\,x\in\Bbb R^2\setminus\Sigma.
$$

Given $d\in S^1$ choose $\vartheta\in S^1$ in such a way that $\vartheta^{\perp}=d$.
Let $x_0\in\Sigma$ and $w=w_{\Sigma}(x;-d,k,x_0)$
be the unique solution of the scattering problem:
$$\begin{array}{c}
\displaystyle
(\triangle+k^2)w=0\,\,\text{in}\,\Bbb R^2\setminus\Sigma,\\
\\
\displaystyle
\frac{\partial w^{\pm}}{\partial\nu^{\pm}}=-\frac{\partial}{\partial\nu^{\pm}}\{(x_0-x)\cdot\vartheta\,e^{-ikx\cdot d}\}\,\,\text{on}\,
\,\Sigma,\\
\\
\displaystyle
\lim_{r\longrightarrow\infty}\,\sqrt{r}\left(\frac{\partial w}{\partial r}-ik\,w\right)=0.
\end{array}
$$
Define
$$
\displaystyle
u_{\Sigma}(x;-d,k,x_0)=(x_0-x)\cdot\vartheta\,e^{-ikx\cdot d}+w_{\Sigma}(x;-d,k,x_0).
$$

Let $R_1>R$ and $y\in\partial B_{R_1}$.
Define
$$\displaystyle
J_{\Sigma}(\tau;\omega,y,k)
=\int_{\partial B_{R}}
\left(\frac{\partial}{\partial\nu}\Phi_{\Sigma}(x,y)\cdot v_{\tau}(x;\omega)-
\frac{\partial}{\partial\nu} v_{\tau}(x;\omega)\cdot\Phi_{\Sigma}(x,y)\right)dS.
$$

The following theorem is what we call an extension of Theorem 1.3 to thin obstacles.

\proclaim{\noindent Theorem 3.2.} Let $\omega$ be regular with respect to $\Sigma$
and let $x_0\in\Sigma$ be the point with $x_0\cdot\omega=h_{\Sigma}(\omega)$.
Assume that $y\in\partial B_{R_1}$ satisfies $u_{\Sigma}(y;-d,k,x_0)\not=0$
for a direction $d\in S^1$ that meets at $x_0$ along a $\Sigma_j$.
Then the formula
$$\displaystyle
\lim_{\tau\longrightarrow\infty}\frac{1}{\tau}\log\vert J_{\Sigma}(\tau;\omega,y,k)\vert
=h_{\Sigma}(\omega),
$$
is valid. Moreover, we have:

if $t\ge h_{\Sigma}(\omega)$, then $\lim_{\tau\longrightarrow\infty}e^{-\tau t}\vert J_{\Sigma}(\tau;\omega,y,k)\vert=0$;

if $t<h_{\Sigma}(\omega)$, then $\lim_{\tau\longrightarrow\infty}e^{-\tau t}\vert J_{\Sigma}(\tau;\omega,y,k)\vert=\infty$.

\endproclaim

The proof of Theorem 3.2 is based on the convergent series expansion of $\Phi_{\Sigma}(\,\cdot\,,y)$ at a corner
or end point of $\Sigma$.  See Propositions 4.4 and 4.5 in \cite{IE4}.
Since the proof of Theorem 3.2 can be done along the same line with that of Theorem 1.3,
we omit the description.

\subsection{Obstacle in a layered medium}

We consider a medium that consists of two parts. One is given by $\Bbb R^2\setminus\overline B_R$
and another is $B_R$.  We assume that the propagation speeds of wave in two parts can be
different from each other.  An obstacle $D$ is embedded in $B_R$ as before.
We assume that $D$ is polygonal.
Let us describe a mathematical formulation of the problem.

Define
$$\gamma(x)=
\left\{\begin{array}{lr}
\displaystyle
\gamma_+, & \quad\text{$x\in\Bbb R^2\setminus B_R$,}\\
\\
\displaystyle
\gamma_-, & \quad\text{$x\in\,B_R$,}
\end{array}
\right.
$$
where $\gamma_{\pm}$ are {\it known} positive constants.

Fix $y\in\Bbb R^2\setminus\overline B_R$.
Set $k_+=k/\sqrt{\gamma_+}$.
Define
$$\displaystyle
\Phi_{\gamma}(x,y)=\Phi_+(x,y)+\epsilon_{\gamma}(x,y),\,\,x\in\Bbb R^2
$$
where
$$\displaystyle
\Phi_+(x,y)=\frac{i}{4}H^{(1)}_0(k_+\vert x-y\vert)
$$
and $\epsilon=\epsilon_{\gamma}$ solves
$$\begin{array}{c}
\displaystyle
(\nabla\cdot\gamma\nabla+k^2)\epsilon
=-\nabla\cdot(\gamma-\gamma_+)\nabla\Phi_+(x,y)\,\,\text{in}\,\Bbb R^2,\\
\\
\displaystyle
\lim_{r\longrightarrow\infty}
\sqrt{r}
\left(\frac{\partial \epsilon}{\partial r}-ik_+\epsilon\right)=0.
\end{array}
$$

Let $E=E_{D,\gamma}(x,y),\,x\in\Bbb R^2\setminus\overline D$ solve
$$\begin{array}{c}
\displaystyle
(\nabla\cdot\gamma\nabla +k^2)E=0\,\,\text{in}\,\Bbb R^2\setminus\overline D,\\
\\
\displaystyle
\gamma\frac{\partial E}{\partial\nu}=-\gamma\frac{\partial}{\partial\nu}\Phi_{\gamma}(x,y)\,\,\text{on}\,\partial D,\\
\\
\displaystyle
\lim_{r\longrightarrow\infty}
\sqrt{r}
\left(\frac{\partial E}{\partial r}-ik_+E\right)=0.
\end{array}
$$
Define
$$\displaystyle
\Phi_{D,\gamma}(x,y)=\Phi_{\gamma}(x,y)+E_{D,\gamma}(x,y).
$$

Note that the existence and uniqueness of the solutions $\epsilon$ and $E$ can be established
by using a variational formulation, for example, see \cite{H}.

Define
$$\displaystyle
K(\tau; \omega, y, k)
=
\int_{\partial B_{R}}
\left(
\gamma_{-}\frac{\partial}{\partial\nu}\Phi_{D,\gamma}^-(x,y)\cdot v_{\tau}^-(x;\omega)
-\gamma_{-}\frac{\partial}{\partial\nu}v_{\tau}^-(x;\omega)\cdot\Phi_{D,\gamma}^-(x,y)\right)dS(x),
\tag {3.1}
$$
where
$$\begin{array}{c}
\displaystyle
v_{\tau}^-(x;\omega)=e^{x\cdot(\tau\omega+i\sqrt{\tau^2+k^2_{-}}\omega^{\perp})},\\
\\
\displaystyle
\Phi_{D,\gamma}^-(x,y)=\Phi_{D,\gamma}(x,y),\,\,x\in \overline B_R\setminus\overline D
\end{array}
$$
and $k_{-}=k/\sqrt{\gamma_{-}}$.

\proclaim{\noindent Theorem 3.3.}
Assume that $\omega$ is regular with respect to $D$ and that
$$\displaystyle
\text{diam}\,D<\text{dist}\,(D,\partial B_{R}).
\tag {3.2}
$$
Moreover assume that there exists a $j\in\{1,\cdots,m\}$ such that $k_{-}^2$ is not a Neumann eigenvalue
for $-\triangle $ in $D_j$.
It holds that
$$\displaystyle
\lim_{\tau\longrightarrow\infty}
\frac{1}{\tau}\log
\left\vert K(\tau; \omega, y, k)
\right\vert
=h_D(\omega).
$$
Moreover, we have the following:

if $t\ge h_D(\omega)$, then
$\displaystyle
\lim_{\tau\longrightarrow\infty}e^{-\tau t}\left\vert
K(\tau; \omega, y, k)
\right\vert=0$;

if $t<h_D(\omega)$, then
$\displaystyle
\lim_{\tau\longrightarrow\infty}e^{-\tau t}\left\vert
K(\tau; \omega, y, k)
\right\vert=\infty$.

\endproclaim

{\it\noindent Proof.} Instead of $u(x)=\Phi_D(x,y)$ in the prof of
Theorem 1.2 set $u(x)=\Phi_{D,\gamma}(x,y)$. For this $u$ the same
argument with (3.2) instead of (1.1) as described in Cases A and B
in the proof of Theorem 1.2 works. Thus we have a continuation
$\tilde{u}$ of $u$ in $B_R\setminus\overline D$ onto $B_R$ as a
solution of the Helmholtz equation
$\triangle\tilde{u}+k_{-}^2\tilde{u}=0$ in $B_R$.  Since $u$
satisfies the Neumann boundary condition $\partial
u/\partial\nu=0$ on $\partial D$, $\tilde{u}$ also satisfies the
condition on $\partial D$ and by the assumption $k_{-}^2$, it must
hold that $\tilde{u}=0$ in $D_j$ for some $j$.  Then, the unique
continuation theorem for the Helmholtz equation yields $u=0$ in
$B_R\setminus\overline D$.  This yields that the Cauchy data of
$\Phi_{D,\gamma}(x,y)$ on $\partial B_R$ vanish and thus
$\Phi_{D,\gamma}(x,y)=0$ for $x\in\,\Bbb R^2\setminus(\overline
B_R\cup\{y\})$ by the uniqueness of the Cauchy problem for the
Helmholtz equation with wave number $k^+$.  Since
$\Phi_{D,\gamma}(x,y)$ is singular as $x\longrightarrow y$, this
is a contradiction.

\noindent
$\Box$

In this theorem the data are given by the Cauchy data of the total
wave field on $\partial B_R$. It is an interesting open problem
when the receivers are located on $\partial B_{R+\epsilon}$ with a
$\epsilon>0$ how one can apply the enclosure method in an {\it
explicit form}.

Here we propose one {\it heuristic} approach based on Theorem 3.3 in the case when $\epsilon$ is sufficiently small.

Assume that we have $\Phi_{D,\gamma}(x,y)$ for all $x\in\partial B_{R+\epsilon}$ exactly.
Solve the exterior problem in $\Bbb R^2\setminus\overline B_{R+\epsilon}$:
$$\begin{array}{c}
\displaystyle
\triangle \Psi+k_+^2\Psi=0\,\,\text{in}\,\Bbb R^2\setminus\overline B_{R+\epsilon},\\
\\
\displaystyle
\Psi(x)=\Phi_{D,\gamma}(x,y)-\Phi^+(x,y)\,\,\text{on}\,\partial B_{R+\epsilon},\\
\\
\displaystyle
\lim_{r\longrightarrow\infty}\sqrt{r}\left(\frac{\partial\Psi}{\partial r}-ik_+\Psi\right)=0.
\end{array}
$$
Then we have $\Psi(x)=\Phi_{D,\gamma}(x,y)-\Phi^+(x,y)$ for $x\in\Bbb R^2\setminus B_{R+\epsilon}$.
This gives
$$
\displaystyle
\frac{\partial\Phi_{D,\gamma}}{\partial\nu}(x+\epsilon\nu(x),y)=\frac{\partial\Phi^+}{\partial\nu}(x+\epsilon\nu(x),y)
+\frac{\partial\Psi}{\partial\nu}(x+\epsilon\nu(x)),\,\,x\in\partial B_R.
$$
We use for the computation of the Cauchy data of $\Phi_{D,\gamma}(x,y)$ on $\partial B_R$
from {\it outside} $B_R$ the approximation:
$$\begin{array}{c}
\Phi_{D,\gamma}(x,y)\approx \Phi_{D,\gamma}(x+\epsilon\nu(x),y),\\
\\
\displaystyle
\frac{\partial\Phi_{D,\gamma}}{\partial\nu}(x,y)\approx
\frac{\partial\Phi^+}{\partial\nu}(x+\epsilon\nu(x),y)
+\frac{\partial\Psi}{\partial\nu}(x+\epsilon\nu(x)).
\end{array}
$$
Using these computed Cauchy data from outside $B_R$ and the transmission condition
$$
\displaystyle
\Phi_{D,\gamma}(x,y)=\Phi_{D,\gamma}^-(x,y),\,\,
\gamma_+\frac{\partial\Phi_{D,\gamma}}{\partial\nu}(x,y)
=\gamma_{-}\frac{\partial\Phi_{D,\gamma}^-}{\partial\nu}(x,y),\,x\in\partial B_R
$$
which is implicitly included in the governing equation, we compute
$K(\tau;\omega,y,k)$ by replacing $\Phi_{D,\gamma}^-(x,y)$ and
$\gamma_-(\partial\Phi_{D,\gamma}^-/\partial\nu)(x,y)$ in the
right-hand side of (3.1) with
$\Phi_{D,\gamma}(x+\epsilon\nu(x),y)$ and
$\gamma_+\{(\partial\Phi^+/\partial\nu)(x+\epsilon\nu(x),y)
+(\partial\Psi/\partial\nu)(x+\epsilon\nu(x))\}$, respectively.
Clearly, the effective range of $\tau$ shall depend on the size of
$\epsilon$.

It would be interesting to test this approach numerically and check its performance.
This belongs to a next research plan.

\subsection{Unsolvability of the far-field equation for polygonal obstacles}

Let $k>0$ and $d\in S^1$.  Let $F_D(\varphi;d,k)$ denote the far-field pattern of the scattered
wave $w(x)=u(x;d,k)-e^{ikx\cdot d}$.

Given $y\in\Bbb R^2$ the {\it far-field equation} for unknown $g\in L^2(S^1)$
$$\displaystyle
\int_{S^1}F_D(\varphi;d,k)g(d)dS(d)
=\frac{e^{i\pi/4}}{\sqrt{8\pi k}}e^{-ik\varphi\cdot y},\,\,\varphi\in S^1
\tag {3.3}
$$
plays the central role in the {\it linear sampling method} \cite{CK1}.

Note that the right-hand side of (3.3) coincides with the {\it far-field pattern}
of the field $\Phi_0(x,y)$ with $x=r\varphi$ as $r\longrightarrow\infty$;
the left-hand side of (3.3) coincides with the far-field pattern
of the scattered field $w=w_g$ which is the unique solution
of the scattering problem:
$$\begin{array}{c}
\displaystyle
(\triangle+k^2)w=0\,\,\text{in}\,\Bbb R^2\setminus\overline D,\\
\\
\displaystyle
\frac{\partial w}{\partial\nu}=-\frac{\partial v_g}{\partial\nu}\,\,\text{on}\,\partial D,\\
\\
\displaystyle
\lim_{r\longrightarrow\infty}\sqrt{r}\left(\frac{\partial w}{\partial r}-ikw\right)=0,
\end{array}
$$
where $v_g$ denotes the {\it Herglotz wave function} with density $g$:
$$\displaystyle
v_g(x)=\int_{S^1} e^{ikx\cdot d}g(d)dS(d),\,x\in\Bbb R^2.
$$
Note that $w_g$ satisfies $w_g\vert_{B_R}\in H^1(B_R\setminus\overline D)$ for a sufficiently large $R$
and the inhomogeneous Nuemann boundary condition on $\partial D$ should be considered in a weak sense.

In this section, using the idea of the proof of Theorem 1.2, we give a proof of unsolvability of
equation (3.3) for any $k>0$ and $y\in\Bbb R^2$ provided $D$ is {\it polygonal}.

\proclaim{\noindent Theorem 3.4.}
For any $k>0$ and $y\in\Bbb R^2$ there exists
no solution $g$ of equation (3.3).
\endproclaim

{\it\noindent Proof.}
We employ a contradiction argument.  Assume that equation (3.3) admits a solution $g$.
Then the coincidence of both fa-field patterns of $w_g$ and $\Phi_0(\,\cdot\,,x)$ yields
$w_g(x)=\Phi_0(x,y)$ for $x\in\Bbb R^2\setminus B_R$ with a sufficiently large $R$.
From the unique continuation property for the Helmholtz equation this coincidence gives
$$\displaystyle
w_g(x)=\Phi_0(x,y),\,\,\forall x\in(\Bbb R^2\setminus\overline D)\setminus\{y\}.
\tag {3.4}
$$
Since $w_g\vert_{B_R}\in H^1(B_R\setminus\overline D)$ and
$\Phi_0(\,\cdot\,, y')\vert_{B_R}$ does not belong to
$H^1(B_R\setminus\overline D)$ for all $y'\in\Bbb R^2\setminus D$
from (3.4) one gets $y\in D$. Note that this part or this type of
argument is well known in the linear sampling method. It shows
that if (3.3) is solvable, then $y\in D$. The problem is the next to
the intermediate conclusion $y\in D$. Now we have
$$\displaystyle
w_g(x)=\Phi_0(x,y),\,\,\forall x\in\,\Bbb R^2\setminus\overline D.
\tag {3.5}
$$
Define
$$\displaystyle
u_g(x)=v_g(x)+w_g(x),\,\,x\in\Bbb R^2\setminus\overline D.
$$
Note that $u=u_g$ satisfies the Helmholtz equation in $\Bbb R^2\setminus\overline D$
and the homogeneous
Neumann boundary condition $\partial u/\partial\nu=0$ on $\partial D$.
From (3.5) we have
$\displaystyle
u_g(x)=v_g(x)+\Phi_0(x,y),\,\,x\in\Bbb R^2\setminus\overline D$ and
this right-hand side gives a continuation $\tilde{u_g}$ of $u_g$ onto $\Bbb R^2\setminus\{y\}$
as a solution of the Helmholtz equation.

Choose a $\omega\in S^1$ that is regular with respect to $D$
and define
$$\displaystyle
I(\tau)=
\int_{\partial B_R}
\left(\frac{\partial u_g}{\partial\nu} v_{\tau}-\frac{\partial v_{\tau}}{\partial\nu}u_g\right)dS,\,\,\tau>0,
\tag {3.6}
$$
where $v_{\tau}(x)=e^{x\cdot(\tau\omega+i\sqrt{\tau^2+k^2}\omega^{\perp})}$.

Let $x_0\in\partial D$ with $x_0\cdot\omega=h_D(\omega)$. One may
assume that $y\in D_1$, where $D_1$ is a connected component of $D$.
Then one can choose a small $\delta>0$
such that if $\vert x-y\vert\le\delta$, then $x\in D_1$ and
$x\cdot\omega<h_D(\omega)-\delta$. Replacing $u_g$ in (3.6) with
$\tilde{u_g}$ and applying integration by parts, we obtain, as
$\tau\longrightarrow\infty$
$$\begin{array}{c}
\displaystyle
e^{-\tau h_D(\omega)}I(\tau)=
e^{-\tau h_D(\omega)}\int_{\partial D}
\left(\frac{\partial\tilde{u_g}}{\partial\nu} v_{\tau}-\frac{\partial v_{\tau}}{\partial\nu}\tilde{u_g}\right)dS\\
\\
\displaystyle
=e^{-\tau h_D(\omega)}\int_{\vert x-y\vert=\delta}
\left(\frac{\partial\tilde{u_g}}{\partial\nu} v_{\tau}-\frac{\partial v_{\tau}}{\partial\nu}\tilde{u_g}\right)dS
=O(\tau e^{-\tau\delta}).
\end{array}
\tag {3.7}
$$
Hereafter we make use of the same notation as those of the proof of Theorem
1.2.

Recalling boundary condition $\partial u_g/\partial\nu=0$ on $\partial D$,
one has the expansion
$$\displaystyle
u_g(r,\theta)=\alpha_1J_0(kr)+\sum_{n=2}^{\infty}\alpha_n J_{\lambda_n}(kr)
\cos\,\lambda_n\theta,\,\,0<r<\eta,\,\,0<\theta<\Theta
$$
and applying the argument for deriving (2.1), we obtain
$$\displaystyle
I(\tau)
\,e^{-i\sqrt{\tau^2+k^2}x_0\cdot\omega^{\perp}}e^{-\tau h_D(\omega)}
\sim
-i\sum_{n=2}^{\infty}\frac{e^{i\frac{\pi}{2}\lambda_n}k^{\lambda_n}\alpha_nK_n}
{s^{\lambda_n}},
\tag {3.8}
$$
where $K_n$ are constants given by the formula $\displaystyle
K_n=e^{ip\lambda_n}+(-1)^n e^{iq\lambda_n}$ and $s=\sqrt{\tau^2+k^2}+\tau$.
Since (3.7) implies that $e^{-\tau h_D(\omega)}I(\tau)$ is rapidly decreasing
as $\tau\longrightarrow\infty$,
all the coefficients of the right-hand side of (3.8) have to vanish, that is
$$\displaystyle
\alpha_n K_n=0,\,\,\forall n\ge 2.
\tag {3.9}
$$

First consider the case when $\Theta/\pi$ is {\it irrational}.
It is easy to see that $K_n\not=0$ for all $n\ge 2$.  Thus from (3.9) one gets
$\alpha_n=0$ and this yields $u_g(r,\theta)=\alpha_1J_0(kr)$
for $0<r<\eta$ and $0<\theta<\Theta$.
Since this right-hand side is an entire solution of the
Helmholtz equation, the unique continuation property of the
solution of the Helmholtz equation yields $\tilde{u_g}(x)=\alpha_1J_0(k\vert
x-x_0\vert)$ in $\Bbb R^2\setminus\overline\{y\}$ and thus one gets
$$\displaystyle
\Phi_0(x,y)=\alpha_1J_0(k\vert x-x_0\vert)-v_g(x),\,\,x\not=y.
$$
Comparing the behaviour as $x\longrightarrow y$ on both sides,
we obtain a contradiction.

Next consider the case when $\Theta/\pi$ is {\it rational}.
Applying the same argument for the derivation of (2.4), we have
a continuation of $u_g$ onto $(\Bbb R^2\setminus\overline D)\cup B_{\eta}(x_0)$
as a solution of the Helmholtz equation and its continuation which we denote by $\tilde{u}'$
satisfies the rotation invariance
$$\displaystyle
\tilde{u}'\left(r,\theta+\frac{2\pi}{a}\right)
=\tilde{u}'(r,\theta), 0<r<\eta,\,\,\theta\in\Bbb R,
$$
where $a\ge 2$ is an integer.
Since the unique continuation property gives $\tilde{u_g}(r,\theta)=\tilde{u}'(r,\theta)$ for $0<r<\eta$ and thus
one gets
$$\displaystyle
\tilde{u_g}\left(r,\theta+\frac{2\pi}{a}\right)
=\tilde{u_g}(r,\theta), 0<r<\eta,\,\theta\in\Bbb R.
\tag {3.10}
$$
Since $\tilde{u_g}$ satisfies the Helmholtz equation for $\vert x-x_0\vert<\vert x_0-y\vert$,
it follows from the unique continuation property and the rotation invariance of the Helmholtz equation
that $\eta$ in (3.10) can be replaced with $\vert x_0-y\vert$:
$$\displaystyle
\tilde{u_g}(x_0+rz(\theta))
=\tilde{u_g}\left(x_0+rz\left(\theta+\frac{2\pi}{a}\right)\right),\,\,0<r<\vert x_0-y\vert,\,\,\theta\in\Bbb R,
\tag {3.11}
$$
where
$\displaystyle
z(\theta)=\cos\,\theta\,\mbox{\boldmath $a$}+\sin\,\theta\,\mbox{\boldmath $a$}^{\perp}$.
Now choose a
$\theta_0$ in such a way that $y=x_0+\vert
y-x_0\vert\,z(\theta_0)$. Since $2\pi/a\le\pi$, we have
$y\not=x_0+\vert y-x_0\vert\,z(\theta_0+2\pi/a)$.  Then letting
$\theta=\theta_0$ and $r\uparrow \vert y-x_0\vert$ in (3.11), we
have a contradiction since $\tilde{u_g}(x)=v_g(x)+\Phi_0(x,y)$ for $x\not=y$ and
$$\displaystyle
\Phi_0(x,y)\sim\frac{1}{2\pi}\log\frac{1}{\vert x-y\vert}
$$
as $x\longrightarrow y$.

\noindent
$\Box$

Using a variational formulation in, e.g., \cite{H}, one can formulate
and establish the unique solvability of the scattering problem of
acoustic wave by a sound-hard obstacle $D$ with Lipschitz boundary.
We use the same notation as those in the case when $D$ is polygonal.
Having the far-field pattern for $D$ with Lipschitz boundary,
one can extend Theorem 3.4 to a slightly general case.
We say that $D$ with a Lipschitz boundary has a {\it horn}, if there exist a $\omega\in S^1$
that is regular with respect to $D$ and $\delta>0$ such that the set
$V\equiv\{x\in D\,\vert\,x\cdot\omega>h_D(\omega)-\delta\}$
becomes a finite cone with the vertex at the point in $\{x\,\vert\,x\cdot\omega=h_D(\omega)\}\cap\partial D$
and the base on $x\cdot\omega=h_D(\omega)-\delta$.
We say that $V$ is a horn.

The conclusion of this section is the following statement and since the proof is really a minor modification
of that of Theorem 3.4 we omit the description of the proof.

\proclaim{\noindent Corollary 3.1.}
If $D$ with a Lipschitz boundary has a horn, then for any $k>0$ and $y\in\Bbb R^2$ there exists
no solution $g$ of equation (3.3).

\endproclaim

Note that, in \cite{LP} the far-field equation for a single
circular obstacle with an arbitrary radius has been considered and
it is shown that the equation is not solvable except for its
center point. Corollary 3.1 means that the existence of a horn $V$
even it is {\it small} prevents the existence of solution of
(3.3) for any $k>0$ and $y\in\Bbb R^2$.

It should be pointed out that the linear sampling method is not based on the solvability
of the far-field equation.  Instead a family of approximate solutions of the far-field equation 
is taken.  See \cite{AT} for interesting study of the method itself.

\centerline{{\bf Acknowledgements}}

This research was partially supported by the Grant-in-Aid for
Scientific Research (C)(No. 21540162) of Japan  Society for the
Promotion of Science.

\vskip1cm
\noindent
e-mail address

ikehata@math.sci.gunma-u.ac.jp

\end{document}